\RequirePackage{filecontents}

\documentclass[a4paper, 11pt, reqno]{amsart}
\usepackage{ifxetex}
\ifxetex
    \PassOptionsToPackage{no-math}{fontspec} 
\fi
\usepackage{kotex}
\usepackage{a4wide}
\usepackage{amssymb}
\usepackage{tikz}
\usepackage{graphicx}

\theoremstyle{plain}
\newtheorem{thm}{Theorem}

\newtheorem{lem}[thm]{Lemma}
\newtheorem{prop}[thm]{Proposition}

\theoremstyle{definition}

\theoremstyle{remark}

\newtheorem*{eg}{Example}

\newcommand\set[1]{\left\{#1\right\}}
\newcommand\abs[1]{\left|#1\right|}

\begin{document}
\title{The chromatic polynomial for cycle graphs}
\author[J. Lee]{Jonghyeon Lee} 
\address[Jonghyeon Lee]{
Department of Mathematics
\\Inha University
\\Incheon 22212, Korea}
\email{orie73@naver.com}

\author[H. Shin]{Heesung Shin$^\dagger$} 
\address[Heesung Shin]{
Department of Mathematics
\\Inha University
\\Incheon 22212, Korea}
\email{shin@inha.ac.kr}
\date{\today}
\thanks{$\dagger$ Corresponding author. This work was supported by the National Research Foundation of Korea(NRF) grant funded by the Korea government(MSIP) (No. 2017R1C1B2008269).}

\begin{abstract}
Let $P(G,\lambda)$ denote the number of proper vertex colorings of $G$ with $\lambda$ colors.
The chromatic polynomial $P(C_n,\lambda)$ for the cycle graph $C_n$ is well-known as 
$$P(C_n,\lambda) = (\lambda-1)^n+(-1)^n(\lambda-1)$$
for all positive integers $n\ge 1$.
Also its inductive proof is widely well-known by the \emph{deletion-contraction recurrence}.
In this paper, we give this inductive proof again and three other proofs of this formula of the chromatic polynomial for the cycle graph $C_n$.
\end{abstract}
\maketitle

\section{Introduction}
\label{sec:intro}
The number of proper colorings of a graph with finite colors was introduced only for planar graphs by George David Birkhoff \cite{Bir12} in 1912, in an attempt to prove the four color theorem, where the formula for this number was later called by the chromatic polynomial.
In 1932, Hassler Whitney \cite{Whi32} generalized Birkhoff's formula from the planar graphs to general graphs.
In 1968, Ronald Cedric Read \cite{Rea68} introduced the concept of chromatically equivalent graphs and asked which polynomials are the chromatic polynomials of some graph, that remains open.

\subsection*{Chromatic polynomial}
For a graph $G$, a \emph{coloring} means almost always a \emph{(proper) vertex coloring}, which is a labeling of vertices of $G$ with colors such that no two adjacent vertices have the same colors. Let $P(G,\lambda)$ denote the number of (proper) vertex colorings of $G$ with $\lambda$ colors and $\chi(G)$ the least number $\lambda$ satisfying $P(G, \lambda)>0$, where $P(G, \lambda)$ and $\chi(G)$ are called a \emph{chromatic polynomial} and \emph{chromatic number} of $G$, respectively.

In fact, it is clear that the number of $\lambda$-colorings is a polynomial in $\lambda$ from a deletion-contraction recurrence.
\begin{prop}[Deletion-contraction recurrence]
For a given a graph $G$ and an edge $e$ in $G$, we have
\begin{align}
P(G,\lambda) = P(G-e,\lambda) - P(G/e,\lambda),
\label{eq:rec}
\end{align}
where $G-e$ is a graph obtained by deletion the edge $e$
and $G/e$ is a graph obtained by contraction the edge $e$.
\end{prop}

\begin{figure}[t]
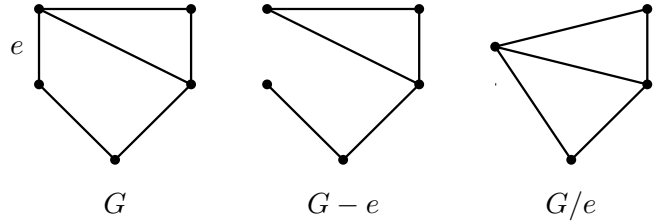

\begin{center}
\begin{pgfpicture}{4.22mm}{-9.86mm}{92.50mm}{22.50mm}
\pgfsetxvec{\pgfpoint{1.00mm}{0mm}}
\pgfsetyvec{\pgfpoint{0mm}{1.00mm}}
\color[rgb]{0,0,0}\pgfsetlinewidth{0.30mm}\pgfsetdash{}{0mm}
\pgfcircle[fill]{\pgfxy(10.00,10.00)}{0.50mm}
\pgfcircle[stroke]{\pgfxy(10.00,10.00)}{0.50mm}
\pgfcircle[fill]{\pgfxy(30.00,10.00)}{0.50mm}
\pgfcircle[stroke]{\pgfxy(30.00,10.00)}{0.50mm}
\pgfcircle[fill]{\pgfxy(30.00,20.00)}{0.50mm}
\pgfcircle[stroke]{\pgfxy(30.00,20.00)}{0.50mm}
\pgfcircle[fill]{\pgfxy(10.00,20.00)}{0.50mm}
\pgfcircle[stroke]{\pgfxy(10.00,20.00)}{0.50mm}
\pgfmoveto{\pgfxy(10.00,10.00)}\pgflineto{\pgfxy(10.00,10.00)}\pgfstroke
\pgfmoveto{\pgfxy(10.00,20.00)}\pgflineto{\pgfxy(10.00,10.00)}\pgfstroke
\pgfmoveto{\pgfxy(10.00,20.00)}\pgflineto{\pgfxy(30.00,20.00)}\pgfstroke
\pgfmoveto{\pgfxy(30.00,20.00)}\pgflineto{\pgfxy(30.00,10.00)}\pgfstroke
\pgfmoveto{\pgfxy(30.00,10.00)}\pgflineto{\pgfxy(10.00,20.00)}\pgfstroke
\pgfmoveto{\pgfxy(20.00,0.00)}\pgflineto{\pgfxy(10.00,10.00)}\pgfstroke
\pgfmoveto{\pgfxy(20.00,0.00)}\pgflineto{\pgfxy(30.00,10.00)}\pgfstroke
\pgfcircle[fill]{\pgfxy(20.00,0.00)}{0.50mm}
\pgfcircle[stroke]{\pgfxy(20.00,0.00)}{0.50mm}
\pgfputat{\pgfxy(8.00,14.00)}{\pgfbox[bottom,left]{\fontsize{11.38}{13.66}\selectfont \makebox[0pt][r]{$e$}}}
\pgfputat{\pgfxy(20.00,-7.00)}{\pgfbox[bottom,left]{\fontsize{11.38}{13.66}\selectfont \makebox[0pt]{$G$}}}
\pgfputat{\pgfxy(80.00,-7.00)}{\pgfbox[bottom,left]{\fontsize{11.38}{13.66}\selectfont \makebox[0pt]{$G/e$}}}
\pgfputat{\pgfxy(50.00,-7.00)}{\pgfbox[bottom,left]{\fontsize{11.38}{13.66}\selectfont \makebox[0pt]{$G-e$}}}
\pgfcircle[fill]{\pgfxy(40.00,10.00)}{0.50mm}
\pgfcircle[stroke]{\pgfxy(40.00,10.00)}{0.50mm}
\pgfcircle[fill]{\pgfxy(60.00,10.00)}{0.50mm}
\pgfcircle[stroke]{\pgfxy(60.00,10.00)}{0.50mm}
\pgfcircle[fill]{\pgfxy(60.00,20.00)}{0.50mm}
\pgfcircle[stroke]{\pgfxy(60.00,20.00)}{0.50mm}
\pgfcircle[fill]{\pgfxy(40.00,20.00)}{0.50mm}
\pgfcircle[stroke]{\pgfxy(40.00,20.00)}{0.50mm}
\pgfmoveto{\pgfxy(40.00,10.00)}\pgflineto{\pgfxy(40.00,10.00)}\pgfstroke
\pgfmoveto{\pgfxy(40.00,20.00)}\pgflineto{\pgfxy(60.00,20.00)}\pgfstroke
\pgfmoveto{\pgfxy(60.00,20.00)}\pgflineto{\pgfxy(60.00,10.00)}\pgfstroke
\pgfmoveto{\pgfxy(60.00,10.00)}\pgflineto{\pgfxy(40.00,20.00)}\pgfstroke
\pgfmoveto{\pgfxy(50.00,0.00)}\pgflineto{\pgfxy(40.00,10.00)}\pgfstroke
\pgfmoveto{\pgfxy(50.00,0.00)}\pgflineto{\pgfxy(60.00,10.00)}\pgfstroke
\pgfcircle[fill]{\pgfxy(50.00,0.00)}{0.50mm}
\pgfcircle[stroke]{\pgfxy(50.00,0.00)}{0.50mm}
\pgfcircle[fill]{\pgfxy(90.00,10.00)}{0.50mm}
\pgfcircle[stroke]{\pgfxy(90.00,10.00)}{0.50mm}
\pgfcircle[fill]{\pgfxy(90.00,20.00)}{0.50mm}
\pgfcircle[stroke]{\pgfxy(90.00,20.00)}{0.50mm}
\pgfcircle[fill]{\pgfxy(70.00,15.00)}{0.50mm}
\pgfcircle[stroke]{\pgfxy(70.00,15.00)}{0.50mm}
\pgfmoveto{\pgfxy(70.00,10.00)}\pgflineto{\pgfxy(70.00,10.00)}\pgfstroke
\pgfmoveto{\pgfxy(70.00,15.00)}\pgflineto{\pgfxy(90.00,20.00)}\pgfstroke
\pgfmoveto{\pgfxy(90.00,20.00)}\pgflineto{\pgfxy(90.00,10.00)}\pgfstroke
\pgfmoveto{\pgfxy(90.00,10.00)}\pgflineto{\pgfxy(70.00,15.00)}\pgfstroke
\pgfmoveto{\pgfxy(80.00,0.00)}\pgflineto{\pgfxy(70.00,15.00)}\pgfstroke
\pgfmoveto{\pgfxy(80.00,0.00)}\pgflineto{\pgfxy(90.00,10.00)}\pgfstroke
\pgfcircle[fill]{\pgfxy(80.00,0.00)}{0.50mm}
\pgfcircle[stroke]{\pgfxy(80.00,0.00)}{0.50mm}
\end{pgfpicture}%
\end{center}
\caption{$G$ , $G-e$ and $G/e$}
\label{fig:rec}
\end{figure}

\begin{eg}
The chromatic polynomials of graphs in Figure~\ref{fig:rec} are
\begin{align*}
P(G,\lambda) &=\lambda(\lambda-1)^2(\lambda-2),\\
P(G-e,\lambda)&=\lambda^2(\lambda-1)(\lambda-2), and\\
P(G/e,\lambda)&=\lambda(\lambda-1)(\lambda-2).
\end{align*}
It is confirmed that \eqref{eq:rec} is true for the graph $G$ and the edge $e$ in Figure~\ref{fig:rec}.
\end{eg}

\subsection*{Cycle graph}
A \emph{cycle graph $C_n$} is a graph that consists of a single cycle of length $n$, which could be drown by a $n$-polygonal graph in a plane.
The chromatic polynomial for cycle graph $C_n$ is well-known as follows.
\begin{thm}
\label{thm:main}
For a positive integer $n\ge1$, the chromatic polynomial for cycle graph $C_n$ is
\begin{align}
P(C_n,\lambda) = (\lambda-1)^n+(-1)^n(\lambda-1)
\label{eq:main}
\end{align}
\end{thm}

\begin{eg}
For an integer $n \le 3$, it is easily checked that the chromatic polynomials of $C_n$ are from \eqref{eq:main} as follows.
\begin{align*}
P(C_1, \lambda) &= (\lambda-1)+(-1)(\lambda-1) = 0 ,\\
P(C_2, \lambda) &= (\lambda-1)^2+(-1)^2(\lambda-1) = \lambda(\lambda-1),\\
P(C_3, \lambda) &= (\lambda-1)^3+(-1)^3(\lambda-1) = \lambda(\lambda-1)(\lambda-2).
\end{align*}

As shown in Figure~\ref{fig:cycle}, the cycle graph $C_1$ is a graph with one vertex and one loop and $C_1$ cannot be colored, that means $P(C_1, \lambda) = 0$.
The cycle graph $C_2$ is a graph with two vertices, where two edges between two vertices,
and $C_2$ can have colorings by assigning two vertices with different colors, that means $P(C_2, \lambda) = \lambda(\lambda-1)$.
The cycle graph $C_3$ is drawn by a triangle and $C_3$ can have colorings by assigning all three vertices with different colors, that means $P(C_3, \lambda) = \lambda(\lambda-1)(\lambda-2)$.

\begin{figure}[t]
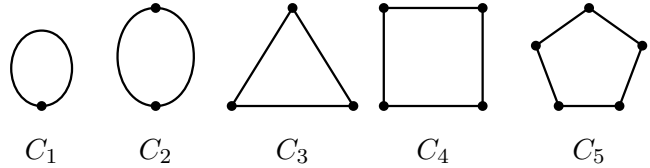

\begin{center}
\begin{pgfpicture}{19.00mm}{-9.86mm}{106.50mm}{15.50mm}
\pgfsetxvec{\pgfpoint{1.00mm}{0mm}}
\pgfsetyvec{\pgfpoint{0mm}{1.00mm}}
\color[rgb]{0,0,0}\pgfsetlinewidth{0.30mm}\pgfsetdash{}{0mm}
\pgfputat{\pgfxy(25.00,-7.00)}{\pgfbox[bottom,left]{\fontsize{11.38}{13.66}\selectfont \makebox[0pt]{$C_1$}}}
\pgfcircle[fill]{\pgfxy(25.00,0.00)}{0.50mm}
\pgfcircle[stroke]{\pgfxy(25.00,0.00)}{0.50mm}
\pgfmoveto{\pgfxy(25.00,0.00)}\pgfcurveto{\pgfxy(22.59,0.00)}{\pgfxy(21.00,2.40)}{\pgfxy(21.00,5.00)}\pgfcurveto{\pgfxy(21.00,7.60)}{\pgfxy(22.59,10.00)}{\pgfxy(25.00,10.00)}\pgfcurveto{\pgfxy(27.41,10.00)}{\pgfxy(29.00,7.60)}{\pgfxy(29.00,5.00)}\pgfcurveto{\pgfxy(29.00,2.40)}{\pgfxy(27.41,0.00)}{\pgfxy(25.00,0.00)}\pgfclosepath\pgfstroke
\pgfputat{\pgfxy(40.00,-7.00)}{\pgfbox[bottom,left]{\fontsize{11.38}{13.66}\selectfont \makebox[0pt]{$C_2$}}}
\pgfcircle[fill]{\pgfxy(40.00,0.00)}{0.50mm}
\pgfcircle[stroke]{\pgfxy(40.00,0.00)}{0.50mm}
\pgfmoveto{\pgfxy(40.00,0.00)}\pgfcurveto{\pgfxy(37.09,0.00)}{\pgfxy(35.20,2.88)}{\pgfxy(35.00,6.00)}\pgfcurveto{\pgfxy(34.77,9.57)}{\pgfxy(36.75,13.00)}{\pgfxy(40.00,13.00)}\pgfcurveto{\pgfxy(43.25,13.00)}{\pgfxy(45.23,9.57)}{\pgfxy(45.00,6.00)}\pgfcurveto{\pgfxy(44.80,2.88)}{\pgfxy(42.91,0.00)}{\pgfxy(40.00,0.00)}\pgfclosepath\pgfstroke
\pgfcircle[fill]{\pgfxy(40.00,13.00)}{0.50mm}
\pgfcircle[stroke]{\pgfxy(40.00,13.00)}{0.50mm}
\pgfmoveto{\pgfxy(70.00,13.00)}\pgflineto{\pgfxy(83.00,13.00)}\pgflineto{\pgfxy(83.00,0.00)}\pgflineto{\pgfxy(70.00,0.00)}\pgfclosepath\pgfstroke
\pgfmoveto{\pgfxy(66.00,0.00)}\pgflineto{\pgfxy(50.00,0.00)}\pgflineto{\pgfxy(58.00,13.00)}\pgfclosepath\pgfstroke
\pgfmoveto{\pgfxy(93.00,0.00)}\pgflineto{\pgfxy(101.00,0.00)}\pgflineto{\pgfxy(104.00,8.00)}\pgflineto{\pgfxy(97.00,13.00)}\pgflineto{\pgfxy(90.00,8.00)}\pgfclosepath\pgfstroke
\pgfcircle[fill]{\pgfxy(93.00,0.00)}{0.50mm}
\pgfcircle[stroke]{\pgfxy(93.00,0.00)}{0.50mm}
\pgfcircle[fill]{\pgfxy(101.00,0.00)}{0.50mm}
\pgfcircle[stroke]{\pgfxy(101.00,0.00)}{0.50mm}
\pgfcircle[fill]{\pgfxy(104.00,8.00)}{0.50mm}
\pgfcircle[stroke]{\pgfxy(104.00,8.00)}{0.50mm}
\pgfcircle[fill]{\pgfxy(97.00,13.00)}{0.50mm}
\pgfcircle[stroke]{\pgfxy(97.00,13.00)}{0.50mm}
\pgfcircle[fill]{\pgfxy(90.00,8.00)}{0.50mm}
\pgfcircle[stroke]{\pgfxy(90.00,8.00)}{0.50mm}
\pgfcircle[fill]{\pgfxy(83.00,0.00)}{0.50mm}
\pgfcircle[stroke]{\pgfxy(83.00,0.00)}{0.50mm}
\pgfcircle[fill]{\pgfxy(83.00,13.00)}{0.50mm}
\pgfcircle[stroke]{\pgfxy(83.00,13.00)}{0.50mm}
\pgfcircle[fill]{\pgfxy(70.00,13.00)}{0.50mm}
\pgfcircle[stroke]{\pgfxy(70.00,13.00)}{0.50mm}
\pgfcircle[fill]{\pgfxy(70.00,0.00)}{0.50mm}
\pgfcircle[stroke]{\pgfxy(70.00,0.00)}{0.50mm}
\pgfcircle[fill]{\pgfxy(66.00,0.00)}{0.50mm}
\pgfcircle[stroke]{\pgfxy(66.00,0.00)}{0.50mm}
\pgfcircle[fill]{\pgfxy(58.00,13.00)}{0.50mm}
\pgfcircle[stroke]{\pgfxy(58.00,13.00)}{0.50mm}
\pgfcircle[fill]{\pgfxy(50.00,0.00)}{0.50mm}
\pgfcircle[stroke]{\pgfxy(50.00,0.00)}{0.50mm}
\pgfputat{\pgfxy(97.00,-7.00)}{\pgfbox[bottom,left]{\fontsize{11.38}{13.66}\selectfont \makebox[0pt]{$C_5$}}}
\pgfputat{\pgfxy(76.50,-7.00)}{\pgfbox[bottom,left]{\fontsize{11.38}{13.66}\selectfont \makebox[0pt]{$C_4$}}}
\pgfputat{\pgfxy(58.00,-7.00)}{\pgfbox[bottom,left]{\fontsize{11.38}{13.66}\selectfont \makebox[0pt]{$C_3$}}}
\end{pgfpicture}%
\end{center}
\caption{$C_n$  $(1 \leq n \leq 5)$}
\label{fig:cycle}
\end{figure}

\end{eg}

%

\section{Four proofs of Theorem~\ref{thm:main}}
In this section, we show the formula~\eqref{eq:main} in four different ways.

\subsection{Inductive proof}
This inductive proof is widely well-known.
A \emph{path graph $P_n$} is a connected graph in which $n-1$ edges connect $n$ vertices of vertex degree at most $2$, which could be drawn on a single straight line.
The chromatic polynomial for path graph $P_n$ is easily obtained by coloring all vertices $v_1, \dots, v_n$ where $v_i$ and $v_{i+1}$ have different colors for $i=1,
\dots, n-1$.
\begin{lem}
\label{lem:path}
For a positive integer $n\ge1$, the chromatic polynomial for path graph $P_n$ is
\begin{align}
P(P_n,\lambda)=\lambda(\lambda-1)^{n-1}.
\label{eq:path}
\end{align}
\end{lem}

We use an induction on the number $n$ of vertices by the deletion-contraction recurrence and the above lemma for path graph:
It is already shown that \eqref{eq:main} is true for $n\le3$ by the example in Section~\ref{sec:intro}.
Assume that \eqref{eq:main} is true for a positive integer $n$.
Using \eqref{eq:rec} and \eqref{eq:path}, we have
\begin{align*}
P(C_{n+1},\lambda)
&=P(C_{n+1}-e,\lambda)-P(C_{n+1}/e,\lambda)
\tag*{by \eqref{eq:rec}} \\
&=P(P_{n+1},\lambda)-P(C_n,\lambda) \\
&=\lambda(\lambda-1)^n-\left( (\lambda-1)^n+(-1)^n(\lambda-1) \right)
\tag*{by \eqref{eq:path}}\\
&=(\lambda-1)^{n+1}+(-1)^{n+1}(\lambda-1).
\end{align*}

\begin{figure}[t]
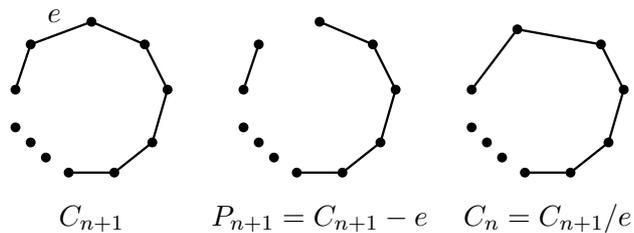

\begin{center}
\centering
\begin{pgfpicture}{7.50mm}{-9.86mm}{92.50mm}{25.14mm}
\pgfsetxvec{\pgfpoint{1.00mm}{0mm}}
\pgfsetyvec{\pgfpoint{0mm}{1.00mm}}
\color[rgb]{0,0,0}\pgfsetlinewidth{0.30mm}\pgfsetdash{}{0mm}
\pgfputat{\pgfxy(16.00,20.00)}{\pgfbox[bottom,left]{\fontsize{11.38}{13.66}\selectfont \makebox[0pt][r]{$e$}}}
\pgfputat{\pgfxy(20.00,-7.00)}{\pgfbox[bottom,left]{\fontsize{11.38}{13.66}\selectfont \makebox[0pt]{$C_{n+1}$}}}
\pgfputat{\pgfxy(80.00,-7.00)}{\pgfbox[bottom,left]{\fontsize{11.38}{13.66}\selectfont \makebox[0pt]{$C_n=C_{n+1}/e$}}}
\pgfputat{\pgfxy(50.00,-7.00)}{\pgfbox[bottom,left]{\fontsize{11.38}{13.66}\selectfont \makebox[0pt]{$P_{n+1}=C_{n+1}-e$}}}
\pgfmoveto{\pgfxy(20.00,20.00)}\pgflineto{\pgfxy(27.00,17.00)}\pgfstroke
\pgfmoveto{\pgfxy(27.00,17.00)}\pgflineto{\pgfxy(30.00,11.00)}\pgfstroke
\pgfmoveto{\pgfxy(30.00,11.00)}\pgflineto{\pgfxy(28.00,4.00)}\pgfstroke
\pgfmoveto{\pgfxy(28.00,4.00)}\pgflineto{\pgfxy(23.00,0.00)}\pgfstroke
\pgfmoveto{\pgfxy(23.00,0.00)}\pgflineto{\pgfxy(17.00,0.00)}\pgfstroke
\pgfmoveto{\pgfxy(20.00,20.00)}\pgflineto{\pgfxy(12.00,17.00)}\pgfstroke
\pgfmoveto{\pgfxy(12.00,17.00)}\pgflineto{\pgfxy(10.00,11.00)}\pgfstroke
\pgfcircle[fill]{\pgfxy(10.00,11.00)}{0.50mm}
\pgfcircle[stroke]{\pgfxy(10.00,11.00)}{0.50mm}
\pgfcircle[fill]{\pgfxy(12.00,4.00)}{0.50mm}
\pgfcircle[stroke]{\pgfxy(12.00,4.00)}{0.50mm}
\pgfcircle[fill]{\pgfxy(14.00,2.00)}{0.50mm}
\pgfcircle[stroke]{\pgfxy(14.00,2.00)}{0.50mm}
\pgfcircle[fill]{\pgfxy(17.00,0.00)}{0.50mm}
\pgfcircle[stroke]{\pgfxy(17.00,0.00)}{0.50mm}
\pgfcircle[fill]{\pgfxy(10.00,6.00)}{0.50mm}
\pgfcircle[stroke]{\pgfxy(10.00,6.00)}{0.50mm}
\pgfcircle[fill]{\pgfxy(23.00,0.00)}{0.50mm}
\pgfcircle[stroke]{\pgfxy(23.00,0.00)}{0.50mm}
\pgfcircle[fill]{\pgfxy(28.00,4.00)}{0.50mm}
\pgfcircle[stroke]{\pgfxy(28.00,4.00)}{0.50mm}
\pgfcircle[fill]{\pgfxy(30.00,11.00)}{0.50mm}
\pgfcircle[stroke]{\pgfxy(30.00,11.00)}{0.50mm}
\pgfcircle[fill]{\pgfxy(27.00,17.00)}{0.50mm}
\pgfcircle[stroke]{\pgfxy(27.00,17.00)}{0.50mm}
\pgfcircle[fill]{\pgfxy(20.00,20.00)}{0.50mm}
\pgfcircle[stroke]{\pgfxy(20.00,20.00)}{0.50mm}
\pgfcircle[fill]{\pgfxy(12.00,17.00)}{0.50mm}
\pgfcircle[stroke]{\pgfxy(12.00,17.00)}{0.50mm}
\pgfmoveto{\pgfxy(50.00,20.00)}\pgflineto{\pgfxy(57.00,17.00)}\pgfstroke
\pgfmoveto{\pgfxy(57.00,17.00)}\pgflineto{\pgfxy(60.00,11.00)}\pgfstroke
\pgfmoveto{\pgfxy(60.00,11.00)}\pgflineto{\pgfxy(58.00,4.00)}\pgfstroke
\pgfmoveto{\pgfxy(58.00,4.00)}\pgflineto{\pgfxy(53.00,0.00)}\pgfstroke
\pgfmoveto{\pgfxy(53.00,0.00)}\pgflineto{\pgfxy(47.00,0.00)}\pgfstroke
\pgfmoveto{\pgfxy(42.00,17.00)}\pgflineto{\pgfxy(40.00,11.00)}\pgfstroke
\pgfcircle[fill]{\pgfxy(40.00,11.00)}{0.50mm}
\pgfcircle[stroke]{\pgfxy(40.00,11.00)}{0.50mm}
\pgfcircle[fill]{\pgfxy(42.00,4.00)}{0.50mm}
\pgfcircle[stroke]{\pgfxy(42.00,4.00)}{0.50mm}
\pgfcircle[fill]{\pgfxy(44.00,2.00)}{0.50mm}
\pgfcircle[stroke]{\pgfxy(44.00,2.00)}{0.50mm}
\pgfcircle[fill]{\pgfxy(47.00,0.00)}{0.50mm}
\pgfcircle[stroke]{\pgfxy(47.00,0.00)}{0.50mm}
\pgfcircle[fill]{\pgfxy(40.00,6.00)}{0.50mm}
\pgfcircle[stroke]{\pgfxy(40.00,6.00)}{0.50mm}
\pgfcircle[fill]{\pgfxy(53.00,0.00)}{0.50mm}
\pgfcircle[stroke]{\pgfxy(53.00,0.00)}{0.50mm}
\pgfcircle[fill]{\pgfxy(58.00,4.00)}{0.50mm}
\pgfcircle[stroke]{\pgfxy(58.00,4.00)}{0.50mm}
\pgfcircle[fill]{\pgfxy(60.00,11.00)}{0.50mm}
\pgfcircle[stroke]{\pgfxy(60.00,11.00)}{0.50mm}
\pgfcircle[fill]{\pgfxy(57.00,17.00)}{0.50mm}
\pgfcircle[stroke]{\pgfxy(57.00,17.00)}{0.50mm}
\pgfcircle[fill]{\pgfxy(50.00,20.00)}{0.50mm}
\pgfcircle[stroke]{\pgfxy(50.00,20.00)}{0.50mm}
\pgfcircle[fill]{\pgfxy(42.00,17.00)}{0.50mm}
\pgfcircle[stroke]{\pgfxy(42.00,17.00)}{0.50mm}
\pgfmoveto{\pgfxy(76.00,19.00)}\pgflineto{\pgfxy(87.00,17.00)}\pgfstroke
\pgfmoveto{\pgfxy(87.00,17.00)}\pgflineto{\pgfxy(90.00,11.00)}\pgfstroke
\pgfmoveto{\pgfxy(90.00,11.00)}\pgflineto{\pgfxy(88.00,4.00)}\pgfstroke
\pgfmoveto{\pgfxy(88.00,4.00)}\pgflineto{\pgfxy(83.00,0.00)}\pgfstroke
\pgfmoveto{\pgfxy(83.00,0.00)}\pgflineto{\pgfxy(77.00,0.00)}\pgfstroke
\pgfmoveto{\pgfxy(76.00,19.00)}\pgflineto{\pgfxy(70.00,11.00)}\pgfstroke
\pgfcircle[fill]{\pgfxy(70.00,11.00)}{0.50mm}
\pgfcircle[stroke]{\pgfxy(70.00,11.00)}{0.50mm}
\pgfcircle[fill]{\pgfxy(72.00,4.00)}{0.50mm}
\pgfcircle[stroke]{\pgfxy(72.00,4.00)}{0.50mm}
\pgfcircle[fill]{\pgfxy(74.00,2.00)}{0.50mm}
\pgfcircle[stroke]{\pgfxy(74.00,2.00)}{0.50mm}
\pgfcircle[fill]{\pgfxy(77.00,0.00)}{0.50mm}
\pgfcircle[stroke]{\pgfxy(77.00,0.00)}{0.50mm}
\pgfcircle[fill]{\pgfxy(70.00,6.00)}{0.50mm}
\pgfcircle[stroke]{\pgfxy(70.00,6.00)}{0.50mm}
\pgfcircle[fill]{\pgfxy(83.00,0.00)}{0.50mm}
\pgfcircle[stroke]{\pgfxy(83.00,0.00)}{0.50mm}
\pgfcircle[fill]{\pgfxy(88.00,4.00)}{0.50mm}
\pgfcircle[stroke]{\pgfxy(88.00,4.00)}{0.50mm}
\pgfcircle[fill]{\pgfxy(90.00,11.00)}{0.50mm}
\pgfcircle[stroke]{\pgfxy(90.00,11.00)}{0.50mm}
\pgfcircle[fill]{\pgfxy(87.00,17.00)}{0.50mm}
\pgfcircle[stroke]{\pgfxy(87.00,17.00)}{0.50mm}
\pgfcircle[fill]{\pgfxy(76.00,19.00)}{0.50mm}
\pgfcircle[stroke]{\pgfxy(76.00,19.00)}{0.50mm}
\end{pgfpicture}%
\end{center}
\caption{$C_{n+1}$ , $P_{n+1}$ and $C_n$}
\end{figure}

Thus, \eqref{eq:main} is true for all positive integers $n\ge1$.

\subsection{Proof by inclusion-exclusion principle}
The \emph{inclusion-exclusion principle} is a technique of counting the size of the union of finite sets.
\begin{prop}[Inclusion-exclusion principle]
Let $A_1, A_2, \dots, A_n$ be subsets of a finite set $U$.
Then number of elements excluding their union is as follows
\begin{align*}
\abs{\bigcap_{i=1}^n {\overline{A_i}}}
&=\sum_{I \subset [n]} (-1)^{\abs{I}} \abs{\bigcap_{i\in I} A_i}\\
&=\abs{U} - \sum_{i=1}^{n} \abs{A_i} + \sum_{i<j} \abs{A_i \cap A_j} - \dots +(-1)^n \abs{A_1 \cap \cdots \cap A_n}
\end{align*}
where $\overline{A}$ is the complement of $A$ in $U$.
\end{prop}


Considering every condition to assign different colors to two adjacent vertices, for each edge $e$, we define a finite sets of arbitrary (including improper) colorings to assign same color to two adjacent vertices by the edge $e$.

Let $A_i$ be a set of colorings such that two vertices $v_i$ and $v_{i+1}$ are of same color, where $v_{n+1}$ is regarded as $v_1$.
Applying the inclusion-exclusion principle, we can write the following
\begin{align*}
P(C_n,\lambda)
&=\vert U \vert
- \sum_{i=1}^{n} \abs{A_i}
+ \sum_{i<j} \abs{A_i \cap A_j}
+ \cdots
+ (-1)^n \abs{A_1 \cap \dots \cap A_n} \\
&= \lambda^n
- \binom{n}{1} \lambda^{n-1}
+ \binom{n}{2} \lambda^{n-2}
+ \cdots
+ (-1)^{n-1} \binom{n}{n-1} \lambda
+ (-1)^n \lambda \\
&= (\lambda-1)^n - (-1)^n +(-1)^n \lambda \\
&= (\lambda-1)^n+(-1)^n(\lambda-1).
\end{align*}
Thus, \eqref{eq:main} is true for all positive integers $n\ge1$.
	
\subsection{Algebric proof}
\begin{figure}[t]
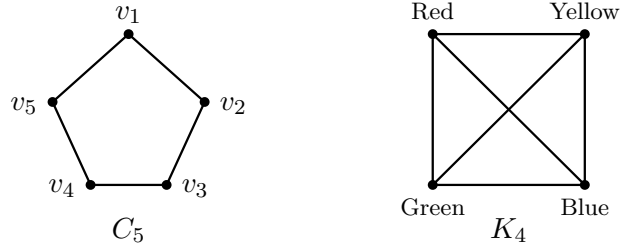

\begin{center}
\centering
\begin{pgfpicture}{2.00mm}{-9.86mm}{86.42mm}{27.14mm}
\pgfsetxvec{\pgfpoint{1.00mm}{0mm}}
\pgfsetyvec{\pgfpoint{0mm}{1.00mm}}
\color[rgb]{0,0,0}\pgfsetlinewidth{0.30mm}\pgfsetdash{}{0mm}
\pgfputat{\pgfxy(20.00,-7.00)}{\pgfbox[bottom,left]{\fontsize{11.38}{13.66}\selectfont \makebox[0pt]{$C_{5}$}}}
\pgfmoveto{\pgfxy(20.00,20.00)}\pgflineto{\pgfxy(10.00,11.00)}\pgflineto{\pgfxy(15.00,0.00)}\pgflineto{\pgfxy(25.00,0.00)}\pgflineto{\pgfxy(30.00,11.00)}\pgfclosepath\pgfstroke
\pgfcircle[fill]{\pgfxy(25.00,0.00)}{0.50mm}
\pgfcircle[stroke]{\pgfxy(25.00,0.00)}{0.50mm}
\pgfcircle[fill]{\pgfxy(30.00,11.00)}{0.50mm}
\pgfcircle[stroke]{\pgfxy(30.00,11.00)}{0.50mm}
\pgfcircle[fill]{\pgfxy(20.00,20.00)}{0.50mm}
\pgfcircle[stroke]{\pgfxy(20.00,20.00)}{0.50mm}
\pgfcircle[fill]{\pgfxy(10.00,11.00)}{0.50mm}
\pgfcircle[stroke]{\pgfxy(10.00,11.00)}{0.50mm}
\pgfcircle[fill]{\pgfxy(15.00,0.00)}{0.50mm}
\pgfcircle[stroke]{\pgfxy(15.00,0.00)}{0.50mm}
\pgfmoveto{\pgfxy(60.00,20.00)}\pgflineto{\pgfxy(60.00,0.00)}\pgflineto{\pgfxy(80.00,0.00)}\pgflineto{\pgfxy(80.00,20.00)}\pgfclosepath\pgfstroke
\pgfmoveto{\pgfxy(80.00,20.00)}\pgflineto{\pgfxy(80.00,20.00)}\pgfstroke
\pgfmoveto{\pgfxy(80.00,20.00)}\pgflineto{\pgfxy(60.00,0.00)}\pgfstroke
\pgfmoveto{\pgfxy(60.00,20.00)}\pgflineto{\pgfxy(80.00,0.00)}\pgfstroke
\pgfcircle[fill]{\pgfxy(80.00,0.00)}{0.50mm}
\pgfcircle[stroke]{\pgfxy(80.00,0.00)}{0.50mm}
\pgfcircle[fill]{\pgfxy(80.00,20.00)}{0.50mm}
\pgfcircle[stroke]{\pgfxy(80.00,20.00)}{0.50mm}
\pgfcircle[fill]{\pgfxy(60.00,20.00)}{0.50mm}
\pgfcircle[stroke]{\pgfxy(60.00,20.00)}{0.50mm}
\pgfcircle[fill]{\pgfxy(60.00,0.00)}{0.50mm}
\pgfcircle[stroke]{\pgfxy(60.00,0.00)}{0.50mm}
\pgfputat{\pgfxy(70.00,-7.00)}{\pgfbox[bottom,left]{\fontsize{11.38}{13.66}\selectfont \makebox[0pt]{$K_{4}$}}}
\pgfputat{\pgfxy(20.00,22.00)}{\pgfbox[bottom,left]{\fontsize{11.38}{13.66}\selectfont \makebox[0pt]{$v_1$}}}
\pgfputat{\pgfxy(32.00,10.00)}{\pgfbox[bottom,left]{\fontsize{11.38}{13.66}\selectfont $v_2$}}
\pgfputat{\pgfxy(27.00,-1.00)}{\pgfbox[bottom,left]{\fontsize{11.38}{13.66}\selectfont $v_3$}}
\pgfputat{\pgfxy(13.00,-1.00)}{\pgfbox[bottom,left]{\fontsize{11.38}{13.66}\selectfont \makebox[0pt][r]{$v_4$}}}
\pgfputat{\pgfxy(8.00,10.00)}{\pgfbox[bottom,left]{\fontsize{11.38}{13.66}\selectfont \makebox[0pt][r]{$v_5$}}}
\pgfputat{\pgfxy(60.00,22.00)}{\pgfbox[bottom,left]{\fontsize{8.54}{10.24}\selectfont \makebox[0pt]{Red}}}
\pgfputat{\pgfxy(80.00,-4.00)}{\pgfbox[bottom,left]{\fontsize{8.54}{10.24}\selectfont \makebox[0pt]{Blue}}}
\pgfputat{\pgfxy(60.00,-4.00)}{\pgfbox[bottom,left]{\fontsize{8.54}{10.24}\selectfont \makebox[0pt]{Green}}}
\pgfputat{\pgfxy(80.00,22.00)}{\pgfbox[bottom,left]{\fontsize{8.54}{10.24}\selectfont \makebox[0pt]{Yellow}}}
\end{pgfpicture}%
\end{center}
\caption{A cycle graph $C_5$ and a graph $K_4$ with names of colors}
\label{fig:one-one}
\end{figure}

Let us consider a case of $n=5$ and $\lambda=4$, that is, to assign the vertices of $C_5$ in four colors: red, blue, yellow, and green. Also let us consider a complete graph $K_4$ with vertex names red, blue, yellow, and green, see Figure~\ref{fig:one-one}.\\

When red-blue-red-yellow-green is assigned in order from the vertex $v_1$ to the vertex $v_5$ in $C_5$, it is corresponding to a closed walk of length $5$ in $K_4$ which begins and ends at red, that is, it is red-blue-red-yellow-green-red in $K_4$.
By generalizing it, we have a correspondence between $\lambda$-colorings of $C_n$ and closed walks of length $n$ in $K_\lambda$. By this correspondence, it is enough to count the number of closed walks of length $n$ in $K_\lambda$, instead of the number of $\lambda$-colorings of $C_n$.

For a graph $G$ with vertex set $\set{v_1, \dots, v_n}$,
the \emph{adjacency matrix} of $G$ is an $n \times n$ square matrix $A$
such that its element $A_{ij}$ is one when there is an edge between two vertices $v_i$ and $v_j$, and zero when there is no edge between $v_i$ and $v_j$.
\begin{figure}[t]
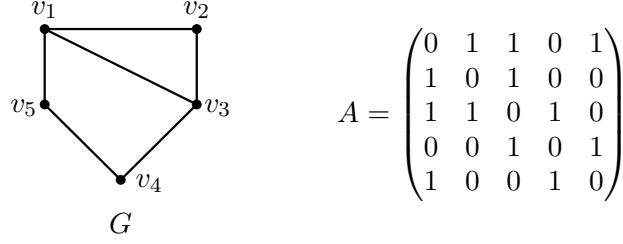

\centering
\begin{minipage}{.3\textwidth}
\centering
\begin{pgfpicture}{3.00mm}{-9.86mm}{37.00mm}{27.14mm}
\pgfsetxvec{\pgfpoint{1.00mm}{0mm}}
\pgfsetyvec{\pgfpoint{0mm}{1.00mm}}
\color[rgb]{0,0,0}\pgfsetlinewidth{0.30mm}\pgfsetdash{}{0mm}
\pgfcircle[fill]{\pgfxy(10.00,10.00)}{0.50mm}
\pgfcircle[stroke]{\pgfxy(10.00,10.00)}{0.50mm}
\pgfcircle[fill]{\pgfxy(30.00,10.00)}{0.50mm}
\pgfcircle[stroke]{\pgfxy(30.00,10.00)}{0.50mm}
\pgfcircle[fill]{\pgfxy(30.00,20.00)}{0.50mm}
\pgfcircle[stroke]{\pgfxy(30.00,20.00)}{0.50mm}
\pgfcircle[fill]{\pgfxy(10.00,20.00)}{0.50mm}
\pgfcircle[stroke]{\pgfxy(10.00,20.00)}{0.50mm}
\pgfmoveto{\pgfxy(10.00,10.00)}\pgflineto{\pgfxy(10.00,10.00)}\pgfstroke
\pgfmoveto{\pgfxy(10.00,20.00)}\pgflineto{\pgfxy(10.00,10.00)}\pgfstroke
\pgfmoveto{\pgfxy(10.00,20.00)}\pgflineto{\pgfxy(30.00,20.00)}\pgfstroke
\pgfmoveto{\pgfxy(30.00,20.00)}\pgflineto{\pgfxy(30.00,10.00)}\pgfstroke
\pgfmoveto{\pgfxy(30.00,10.00)}\pgflineto{\pgfxy(10.00,20.00)}\pgfstroke
\pgfmoveto{\pgfxy(20.00,0.00)}\pgflineto{\pgfxy(10.00,10.00)}\pgfstroke
\pgfmoveto{\pgfxy(20.00,0.00)}\pgflineto{\pgfxy(30.00,10.00)}\pgfstroke
\pgfcircle[fill]{\pgfxy(20.00,0.00)}{0.50mm}
\pgfcircle[stroke]{\pgfxy(20.00,0.00)}{0.50mm}
\pgfputat{\pgfxy(20.00,-7.00)}{\pgfbox[bottom,left]{\fontsize{11.38}{13.66}\selectfont \makebox[0pt]{$G$}}}
\pgfputat{\pgfxy(10.00,22.00)}{\pgfbox[bottom,left]{\fontsize{11.38}{13.66}\selectfont \makebox[0pt]{$v_1$}}}
\pgfputat{\pgfxy(30.00,22.00)}{\pgfbox[bottom,left]{\fontsize{11.38}{13.66}\selectfont \makebox[0pt]{$v_2$}}}
\pgfputat{\pgfxy(31.00,9.00)}{\pgfbox[bottom,left]{\fontsize{11.38}{13.66}\selectfont $v_3$}}
\pgfputat{\pgfxy(22.00,-1.00)}{\pgfbox[bottom,left]{\fontsize{11.38}{13.66}\selectfont $v_4$}}
\pgfputat{\pgfxy(9.00,9.00)}{\pgfbox[bottom,left]{\fontsize{11.38}{13.66}\selectfont \makebox[0pt][r]{$v_5$}}}
\end{pgfpicture}%
\end{minipage}
\quad
\begin{minipage}{.3\textwidth}
$
\displaystyle
A =
\begin{pmatrix}
	0 & 1 & 1 & 0 & 1 \\
	1 & 0 & 1 & 0 & 0 \\
	1 & 1 & 0 & 1 & 0 \\
	0 & 0 & 1 & 0 & 1 \\
	1 & 0 & 0 & 1 & 0
\end{pmatrix}
$
\end{minipage}
\caption{A graph $G$ and its adjacency matrix $A$}
\end{figure}

The following related to an adjacency matrix is well-known.
\begin{prop}
\label{prop:adj}
Let $A$ be the adjacency matrix of the graph $G$ on $n$ vertices $v_1, \dots, v_n$.
Then the $(i,j)$th entry of the matrix $A^n$ is 
the number of the walk of length $n$ beginning at $v_i$ and ending at $v_j$.
\end{prop}

By Proposition~\ref{prop:adj}, we can calculate the number of closed walk of length $n$ in the complete graph $K_\lambda$: Let $A$ be an adjacency matrix of $K_\lambda$. Then $A$ is a $\lambda \times \lambda$ matrix as follows
\begin{align*}
A = \left( a_{ij} \right) =
\begin{pmatrix}
	0 & 1 & \cdots & 1 & 1 \\
	1 & 0 & \cdots & 1 & 1 \\
	\vdots & \vdots & \ddots & \vdots & \vdots \\
	1 & 1 & \cdots & 0 & 1 \\
	1 & 1 & \cdots & 1 & 0
\end{pmatrix},
\end{align*}
where $a_{ij}=0$ if $i=j$, and otherwise $a_{ij}=1$.
So the number of closed walks of length $n$ in $K_\lambda$ is enumerated by $tr(A^n)$,
which equals the sum of all eigenvalues of $A^n$.
Also let all eigenvalues of the matrix $A$ be denoted by $u_1, \dots, u_\lambda$, then all eigenvalues of the matrix $A^n$ are $u_1^n, \dots, u_\lambda^n$.
\begin{align*}
A =
\begin{pmatrix}
	0 & 1 & \cdots & 1 & 1 \\
	1 & 0 & \cdots & 1 & 1 \\
	\vdots & \vdots & \ddots & \vdots & \vdots \\
	1 & 1 & \cdots & 0 & 1 \\
	1 & 1 & \cdots & 1 & 0
\end{pmatrix}
\sim
\begin{pmatrix}
	\lambda-1 & 0 & \cdots & 0 & 0 \\
	0 & -1 & \cdots & 0 & 0 \\
	\vdots & \vdots & \ddots & \vdots & \vdots \\
	0 & 0 & \cdots & -1 & 0 \\
	0 & 0 & \cdots & 0 & -1
\end{pmatrix},
\end{align*}
Since the matrix $A$ have $\lambda$ eigenvalues
$u_1 = \lambda-1$ and $u_2 =\dots=u_\lambda = -1$,
we have
\begin{align*}
tr(A^n) =\sum_{i=1}^{\lambda}u_i^n
= (\lambda-1)^n + \underbrace{(-1)^n + \dots + (-1)^n}_{\text{$\lambda-1$ times}}.
\end{align*}
Thus, \eqref{eq:main} is true for all positive integers $n\ge1$.
			
\subsection{Bijective proof}
Let $X_n$ denote the set of $\lambda$-colorings of $C_n$ and
$[\lambda-1]^n$ be the set of $n$-tuples of positive integers less than $\lambda$,
where $[\lambda-1]$ means $\set{1, \dots, \lambda-1}$.
We consider a mapping $\varphi$ from $\lambda$-colorings of $C_n$ in $X_n$
to $n$-tuples in $[\lambda-1]^n$.

\subsubsection*{A mapping $\varphi$ from $X_n$ to $[\lambda-1]^n$}
The mapping $\varphi:X_n \to [\lambda-1]^n$ is defined as follows:
Let $\omega$ be a $\lambda$-coloring of $C_n$ in $X_n$,
we write $\omega=(\omega_1, \dots, \omega_n)$ where $\omega_i$ is the color of $v_i$ in $C_n$
and it is obvious that $\omega_i \neq \omega_{i+1}$ for $1 \le i \le \lambda$, where $\omega_{n+1}$ is regarded as $\omega_1$.
An entry $\omega_i$ is called a \emph{cyclic descent} of $C$
if $\omega_i > \omega_{i+1}$ for $1\leq i\leq \lambda$.
Then we define $\varphi(\omega) = \sigma = (\sigma_1, \dots, \sigma_n)$ with
\begin{align*}
\sigma_i =
\begin{cases}
\omega_i - 1, & \mbox{\text{if $\omega_i$ is a cyclic descent}}\\
\omega_i, & \mbox{\text{otherwise}}.
\end{cases}
\end{align*}
Given a $\lambda$-coloring $\omega$,
if $\omega_i = \lambda$ then $\omega_{i+1} < \lambda$,
so $\omega_i=\lambda$ should be a cyclic descent.
Thus we have $\sigma_i<\lambda$ for all $1\le i \le n$
and $\varphi(\omega)$ belongs to $[\lambda-1]^n$.

For example, in a case of $n=9$ and $\lambda=4$,
$\omega=(1,2,1,3,2,3,1,4,2) \in X_9$ is given as an example of $4$-colorings of $C_9$.
Here $\omega_2=2$, $\omega_4=3$, $\omega_6 = 3$, $\omega_8=4$, and $\omega_9=2$ are cyclic descents of $\omega$. So we have
$$\varphi(\omega)= \sigma = (1,1,1,2,2,2,1,3,1) \in [3]^9.$$

\subsubsection*{A mapping $\psi$ as the inverse of $\varphi$}
Let $Z_n$ be the set of $n$-tuples $\sigma = (\sigma_1, \sigma_2, \dots, \sigma_n)$ in $[\lambda-1]^n$ with $$\sigma_1 = \sigma_2 = \dots = \sigma_n$$
and it is obvious that the size of $Z_n$ is $\lambda-1$.

We would like to describe a mapping
$\psi: \left( [\lambda-1]^n\setminus Z_n \right) \to X_n$
in order to satisfy $\varphi\circ\psi$ is the identity on $[\lambda-1]^n\setminus Z_n$ as follows:
Given a $\sigma \in [\lambda-1]^n\setminus Z_n$,
we define $\overline{\sigma}=(\overline{\sigma}_1, \dots, \overline{\sigma}_n)$ with
\begin{align*}
\overline{\sigma}_i =
\begin{cases}
\sigma_i + 1, & \mbox{\text{if $\sigma_i$ is a cyclic descent}}\\
\sigma_i, & \mbox{\text{otherwise}}.
\end{cases}
\end{align*}
Since $\overline{\sigma}$ may have consecutive same entries,
we define $\psi(\sigma) = \omega = (\omega_1, \dots, \omega_n)$
from $\overline\sigma$
with $\omega_i = \overline\sigma_i + 1$
for any entry $\overline\sigma_i$ of $\overline\sigma$
with a finite positive even integer $\ell$ satisfying
\begin{align*}
\overline\sigma_{i} = \overline\sigma_{i+1} = \dots = \overline\sigma_{i+\ell-1} \neq \overline\sigma_{i+\ell},
\end{align*}
where $\overline\sigma_{n+k}$ is regarded as $\overline{\sigma}_k$ for $1\le k \le n$,
and $\omega_i = \overline\sigma_i$, otherwise.
Thus $\omega$ has no consecutive same entries
and $1 \le \omega_i \le \lambda$ for all $1\le i \le n$,
so $\psi(\sigma)=\omega$ belongs to $X_n$.
Moreover, it is obvious that $\sigma_i \le \omega_i \le \sigma_i+1$ for all $1\le i \le n$ 
and if $\omega_i = \sigma_i + 1$ for some $1 \le i \le n$ then $\omega_i$ is a cyclic descent in $\omega$.
Hence $\varphi(\omega) = \sigma$ and $\sigma \in [\lambda-1]^n\setminus Z_n$
if and only if $\psi(\sigma)=\omega$.

In a previous example,
$\sigma = (1,1,1,2,2,2,1,3,1)$ is denoted as an example of $9$-tuples in $[3]^9$.
Here $\sigma_6 = 2$, $\sigma_8=3$ are cyclic descents of $\sigma$
and we obtain $\overline\sigma=(1,1,1,2,2,3,1,4,1)$.
And then there exist only three entries $\overline\sigma_2$, $\overline\sigma_4$, and $\overline\sigma_9$ in $\overline\sigma$ satisfying the following
\begin{align*}
k=2:& \quad \overline\sigma_2 = \overline\sigma_3 \neq \overline\sigma_4 \quad (\ell=2),\\
k=4:& \quad \overline\sigma_4 = \overline\sigma_5 \neq \overline\sigma_6 \quad (\ell=2), \text{ and }\\
k=9:& \quad \overline\sigma_9 = \overline\sigma_1 = \overline\sigma_2 = \overline\sigma_3 \neq \overline\sigma_4 \quad (\ell=4),
\end{align*}
so we get $\omega_2 = \overline\sigma_2 + 1 = 2$, $\omega_4= \overline\sigma_4 + 1 =3$, $\omega_9= \overline\sigma_9 + 1 =2$, and
$$\psi(\sigma)= \omega = (1,2,1,3,2,3,1,4,2) \in X_9.$$

Let $Y_n$ be the set of $\lambda$-colorings $\omega$ in $X_n$ with $\varphi(\omega) \in Z_n$.
Since two mapping $\varphi$ and $\psi$ are bijections between $X_n \setminus Y_n$ and $[\lambda-1]^n \setminus Z_n$,
the size of the set $X_n \setminus Y_n$ is same with the size of the $[\lambda-1]^n \setminus Z_n$, which is equal to $(\lambda-1)^n - (\lambda-1).$

When $n$ is even, for any $1\le i \le \lambda-1$,
there exist only two $n$-tuples in $X_n$
\begin{align*}
\omega &= (i+1, i, i+1, i, \dots, i+1, i) \quad \text{ and } \quad
\omega = (i, i+1, i, i+1, \dots, i, i+1)
\end{align*}
satisfying $\varphi(\omega) = (i, i, \dots, i) \in Z_n$.
If $n$ is even, the size of $Y_n$ is equal to $2(\lambda-1)$ and we obtain
\begin{align}
P(C_n,\lambda) &= \abs{X_n}
= \abs{X_n \setminus Y_n} + \abs{Y_n} \tag*{} \\
&= \left[ (\lambda-1)^n - (\lambda-1) \right] + 2(\lambda-1).
\label{eq:even}
\end{align}

When $n$ is odd,
there is no $n$-tuples satisfying $\varphi(\omega) \in Z_n$ and the set $Y_n$ is empty.
If $n$ is odd, we obtain
\begin{align}
P(C_n,\lambda) &= \abs{X_n}
= \abs{X_n \setminus Y_n} + \abs{Y_n} \tag*{} \\
&= \left[ (\lambda-1)^n - (\lambda-1) \right] + 0.
\label{eq:odd}
\end{align}
Therefore, \eqref{eq:main} yields from \eqref{eq:even} and \eqref{eq:odd} for all positive integers $n\ge1$.



\end{document}